\documentclass[12pt]{elsarticle}

\usepackage{amsmath,amssymb,latexsym}
\usepackage{amsfonts}
\usepackage{changepage}
\usepackage[T1]{fontenc} 
\usepackage{textcomp,marvosym}
\usepackage{tabularx}
\usepackage{setspace}
\usepackage{nameref,hyperref}

\usepackage{ulem}
\usepackage[nopatch=eqnum]{microtype}
\usepackage[table]{xcolor}
\usepackage{array}
\usepackage[noabbrev,capitalise]{cleveref}
\usepackage{float}
\usepackage{booktabs}%

\usepackage{url}

\begin{document}

\begin{frontmatter}





\title{A fourth-order dimensionally split exponential time differencing scheme with real and distinct pole rational approximation for non-linear reaction-diffusion systems}


\author[1]{W. K. Attipoe}
\author[2]{E. O. Asante-Asamani \corref{cor1}}

\cortext[cor1]{easantea@clarkson.edu }


\affiliation[1]{organization={College of Western Idaho},
            addressline={5500 E. Opportunity Dr}, 
            city={Nampa},
            postcode={83687}, 
            state={ID},
            country={USA}}
\affiliation[2]{organization={Clarkson University},
            addressline={8 Clarkson Ave}, 
            city={Potsdam},
            postcode={13676}, 
            state={NY},
            country={USA}}

\begin{abstract}
This paper presents a novel fourth-order, L-stable exponential time differencing (ETD) Runge-Kutta type scheme incorporating a dimensional splitting integrating factor method for the numerical solution of nonlinear reaction–diffusion systems. The method, ETDRK4RDP-IF, is obtained by approximating the matrix exponential operators in the classical dimensionally split ETDRK-IF scheme with an accurate fourth-order, L-acceptable rational function with real and distinct poles. The L-acceptability of the rational approximation guarantees effective damping of spurious oscillations arising from nonsmooth initial and/or boundary data and further eliminates complex arithmetic in the solution of the associated linear systems, thereby enhancing computational efficiency and making the scheme well-suited for parallel implementation in low-level programming languages such as Fortran and C++. Numerical experiments confirm the fourth-order accuracy for a range of reaction–diffusion problems with Dirichlet and Neumann boundary conditions, and demonstrate that it outperforms existing standard ETD methods in serial and parallel implementations, yielding CPU time reductions of up to a factor of 20.
\end{abstract}


\end{frontmatter}



\section{Introduction}\label{section1}

Several physical and biological systems that describe the spatio-temporal evolution of chemical substances in a fluid medium can be accurately modeled using time-dependent, nonlinear, multi-dimensional reaction-diffusion equations. Such systems have been studied in medicine \cite{sherratt1990models}, geological engineering \cite{sundararajan2020groundwater}, biological systems, particularly in understanding biological pattern formation \cite{tyson1999model} and monitoring the onset and progression of tumors \cite{anderson2000mathematical} as well as finance \cite{rambeerich2009exponential}. Mathematically they can be represented through the following system of partial differential equations, 

 \begin{eqnarray}\label{eqn1.1}
    \frac{\partial u(x, t)}{\partial t} &=& D \Delta u(x,t) + f(u(x, t), t), \quad \text{in} \ \Omega_T, \ x \in \Omega, \ t \in (0, T], \\ 
     u(x, t) &=& g(x, t), \quad \text{on} \ \Gamma_T,  \nonumber 
 \end{eqnarray}

 \noindent
with suitable initial and boundary data for the domain $\Omega$ and the time interval $[0, T]$. $\Omega$ denotes a bounded open subset of $\mathbb{R}^d$, $d \in \{1,2,3\}$ and $ \Omega_T := \Omega \times (0, T]$ with parabolic boundary $\Gamma_T :=  \overline{\Omega_T} \backslash\Omega_T $. The nonlinear reaction term $f: \Omega_T \longrightarrow \mathbb{R}^s $, is assumed to be a sufficiently smooth function that describes the interactions between all $s$-reacting species. $D = \mathrm{diag}(d_1, d_2, \cdots, d_s), \ d_i > 0 \ \forall i, $ represents the diffusion coefficient, which describes how fast or slow the chemical species diffuse, and the vector-valued solution $u = u(x,t): \overline{\Omega_T} \longrightarrow \mathbb{R}^s$
describes the concentration of each of the $s-$species, which depends on the spatial variable $x = (x_1, x_2, \cdots, x_d) \in \Omega$ and the temporal variable $t \in (0, T]$.

In practice, the initial approach to solving Eq.~\eqref{eqn1.1} involves applying a method of lines technique to discretize the spatial derivatives on a partition of $\Omega$ with $m^d$-nodes to obtain the system of ordinary differential equations 
\begin{eqnarray}\label{eqn1.2}
    \frac{dU}{dt} + AU &=& F(U(t), t), \quad U \in \mathbb{R}^{s\times m^d}, \\ 
     U(x, 0) &=& U_0(x) \nonumber.       
\end{eqnarray}
Here, $A \in \mathbb{R}^{(s\times m^d)\times (s\times m^d)} $ is the discrete approximation of the linear spatial derivative term, $-D\Delta$, and $F \in \mathbb{R}^{(s \times m^d)}$ represents the nonlinear map approximating $f$ on the numerical grid for all $s$-species. $U(t)\in \mathbb{R}^{s\times m^d}$ represents the vector-valued solution across the spatio-temporal grid. 

Numerical solutions to Eq.~\eqref{eqn1.2} pose several numerical challenges, including stiffness induced by spatial diffusion, multiscale dynamics arising from coupled fast and slow reaction processes, and the increased complexity of the nonlinear source term in multidimensional settings $(d\geq 2)$. These issues are compounded by the potential for spurious oscillations resulting from the nonsmooth or incompatible initial and boundary conditions, as well as the need to maintain solution positivity, particularly for biological systems \cite{hundsdorfer2013numerical}.

Exponential time differencing schemes (ETD) are a prominent class of time discretization schemes among a host of others available to solve the resulting stiff ODE systems in Eq.~\eqref{eqn1.2}. Their advantage lies in  accurate treatment of the diffusion term through the matrix exponential, leading to a stable evolution of stiff problems \cite{leveque2007finite}. By applying Dumahel's principle [see, \cite{krabbe1976extension}] to Eq.~\eqref{eqn1.2} and integrating over the interval $[t_0, t]$ we obtain the exact solution as 

\begin{equation}\label{eqn1.3}
    U(t) = \mathrm{e}^{A(t - t_0)}U_0(t) + \int_{t_0}^t \mathrm{e}^{A(t - \tau)} F(U(\tau), \tau) \ \mathrm{d} \tau,
\end{equation}
\noindent
where $t \in [t_0, T]$. Cox and Matthews \cite{cox2002exponential} introduced a class of ETD schemes based on the Runge-Kutta method in 2002 that approximates Eq.~\eqref{eqn1.3} to a higher accuracy. These schemes approximate the integral in Eq.~\eqref{eqn1.3} by interpolating the nonlinear function $F$ with higher-order polynomials and then integrating the resulting expression. In recent developments of the ETDRK schemes, rational functions have been used to approximate the matrix exponential to resolve stability issues with the original Cox and Matthew's formulation as well as create numerical schemes that consist of solving low-cost sparse linear systems of equations (see \cite{asante2025fourth,kleefeld2012etd,yousuf2012numerical,asante2016real,attipoe2025fourth,asante2020second}. These schemes have used both Pad\'e and non-Pad\'e rational functions having real and distinct poles (RDP). The primary advantage of the RDP rational functions is their real and distinct poles which result in schemes devoid of complex arithmetic in solving the resulting linear systems.  


To reduce computational cost and memory requirements, Asante-Asamani, Kleefeld and Wade developed several second-order dimensional splitting techniques, including Strang symmetric splitting, Strang simple splitting\cite{asante2016dimensional}, and an Integrating Factor (IF) \cite{asante2020second} approach for a range of two-dimensional problems having Dirichlet and Neumann boundary data. These methods have been shown to be a competitive alternative to the locally one-dimensional splitting technique proposed by Bhatt and Khaliq \cite{bhatt2015locally} in 2015. Recently, a fourth order version of the integrating factor approach was introduced to improve numerical accuracy\cite{asante2025fourth}. However, the scheme uses an A-acceptable Pad\'e rational function, which is not effective in damping out spurious oscillations at relatively coarse spatio-temporal grids and incurs extra cost in linear solves due to the presence of complex arithmetic. 

In this work, we develop a fourth-order, L-stable, dimensionally splitting ETDRK-IF type scheme by approximating the matrix exponentials in \cite{asante2025fourth} with a fourth-order, non-Pad\'e rational function having real and distinct poles (RDP)\cite{voss1996time}, creating an efficient fourth-order scheme called ETDRK4RDP-IF for the semilinear problem Eq.~\eqref{eqn1.2}. Empirical convergence analysis validates the fourth-order accuracy of the scheme across various systems of RDEs with Dirichlet and Neumann boundary conditions in $\mathbb{R}^2$. The use of RDP rational functions removes complex arithmetic from our linear solves, enabling a more convenient parallelized Fortran implementation of our scheme. The ETDRK4RDP-IF scheme is shown to significantly improve the runtime of the existing fourth-order ETDRK4RDP scheme (without dimensional splitting) and in some systems, achieving up to $5 \times$ and $20 \times$ speed-up in CPU time with serial and parallel implementation, respectively. The performance of the scheme is also compared with ETDRK4P22-IF, which uses a Pad\'e (2,2) rational function to approximate the matrix exponential in the general fourth-order dimensional splitting scheme ETDRK-IF for both serial and parallel implementation.

The rest of the article is organized as follows. In Section(\ref{section2}), we present details on the development of the ETDRK4RDP-IF scheme, along with the RDP rational approximation of the matrix exponential. Serial and parallel algorithms are also discussed for the implementation of the ETDRK4RDP-IF scheme. Numerical experiments showing the performance of the proposed scheme in terms of efficiency and convergence are investigated in Section (\ref{section3}) for various RDE systems. In Section (\ref{section4}), we present the concluding remarks and suggest additional future work.


\section{Development of the ETDRK4RDP-IF scheme}\label{section2}

\subsection{The general fourth-order ETDRK-IF scheme}\label{etdrk4if}
In 2025, Asante-Asamani, Kleefeld and Wade \cite{asante2025fourth} introduced a dimensional splitting version of the class of ETD Runge-Kutta schemes originally developed by Cox-Matthews \cite{cox2002exponential} and Kassam-Trefethen \cite{kassam2005fourth} for solving the semi-linear system of ordinary differential equations in Eq.~\eqref{eqn1.2}. The motivation is to speed up the simulation time when solving problems in high spatial dimensions by decomposing the differentiation matrix into its one-dimensional spatial components. Suppose $A_1$ and $A_2$ are the matrix approximations of $\partial_{xx}$ and $\partial_{yy}$ respectively in 2D, then $A = A_1 + A_2$ is the matrix approximation of the two-dimensional (2D) Laplacian. By writing $B_p$ as the matrix representation of the second-order partial derivative in 1D and $I_p$ a $p-$dimensional identity matrix, then $A_1 = B_p \otimes I_p $ and $A_2 = I_p \otimes B_p $ \cite{hundsdorfer2013numerical}. As established in \cite{asante2020second, asante2016dimensional}, the matrices $A_1$ and $A_2$ commute, which allows for the development of the following scheme



\begin{eqnarray}\label{eqn2.2}
    a_n &=& \mathrm{e}^{-\frac{k}{2}A_2} \mathrm{e}^{-\frac{k}{2}A_1} U_n +\tilde{P}(kA_2) \mathrm{e}^{-\frac{k}{2}A_1}  F(U_n,t_n) \nonumber \\
    b_n &=& \mathrm{e}^{-\frac{k}{2}A_2} \mathrm{e}^{-\frac{k}{2}A_1} U_n +\tilde{P}(kA_2)F \left(a_n,t_n +\frac{k}{2} \right)  \\
    c_n &=& \mathrm{e}^{-\frac{k}{2}A_2} \mathrm{e}^{-\frac{k}{2}A_1} a_n +\tilde{P}(kA_2) \left[ 2 \mathrm{e}^{-\frac{k}{2}A_1} F \left(b_n,t_n+\frac{k}{2} \right)- \mathrm{e}^{-\frac{k}{2}A_1}F(U_n,t_n) \right] \nonumber\\
    U_{n+1} &=& \mathrm{e}^{-kA_1} \mathrm{e}^{-kA_2} U_n + P_1(kA_2) \mathrm{e}^{-kA_1}F(U_n,t_n) \nonumber\\
    && +2P_2(kA_2) \mathrm{e}^{-\frac{k}{2}A_1} G \left( a_n, b_n, t_n + \frac{k}{2} \right) + P_3(kA_2) F(c_n,t_n+k) \nonumber, 
\end{eqnarray}
\noindent
where 
$$\displaystyle G \left( a_n, b_n, t_n + \frac{k}{2} \right) = F\left(a_n, t_n + \frac{k}{2}\right) + F\left(b_n, t_n + \frac{k}{2}\right),$$

\noindent
and 
 \begin{align}\label{eqn2.2T}
    P_1(kA_2)&=  \frac{1}{k^2}(-A_2)^{-3}[-4I+kA_2+\mathrm{e}^{-kA_2}(4I+3kA_2+k^2A_2^2)]\nonumber \\
    P_2(kA_2)&=  \frac{1}{k^2}(-A_2)^{-3}[2I-kA_2-\mathrm{e}^{-kA_2}(2I+kA_2)]\\
    P_3(kA_2)&= \frac{1}{k^2}(-A_2)^{-3}[-4I+3kA_2-k^2A_2^2+\mathrm{e}^{-kA_2}(4I+kA_2)] \nonumber \\
    \tilde{P}(kA_2)&=- A_2^{-1}(\mathrm{e}^{-\frac{k}{2}A_2}-I).\nonumber
\end{align}

Eq.~\eqref{eqn2.2} is the dimensional splitting formulation of the standard fourth-order ETD Runge–Kutta scheme. The system matrix $A,$ appears in terms of the components $A_1$ and $A_2$, which have smaller bandwidth and thus are more computationally efficient to work with.


\subsection{Real and distinct pole rational function (RDP)}\label{rdp}
To approximate the matrix exponentials in the ETDRK-IF scheme Eq.~\eqref{eqn2.2}, we adopt the fourth-order RDP rational approximation to the matrix exponential originally developed by Voss and Khaliq in \cite{voss1996time}. The real-valued rational approximation is expressed as

\begin{equation}\label{eqn2.31}
    R(z) = \frac{N(z)}{D(z)} = \frac{1 + a_1z + a_2z^2 + a_3z^3}{(1-b_1z)(1-b_2z)(1-b_3z)(1-b_4z)}
\end{equation}
where $a_i, b_i \in \mathbb{R}.$ The motivation for requiring the poles and coefficients to be real is to eliminate the need for complex arithmetic within the fully discrete scheme and its associated linear systems, with the expectation of speeding up the run time. Requiring the poles to be distinct leads to the partial fraction decomposition

\begin{equation}\label{eqn2.32}
R(z) = \sum_{i=1}^4 \frac{w_i}{1 - b_iz}. 
\end{equation}

L-acceptability, that is, A-acceptability combined with the property  $|R(z)| \longrightarrow 0$ as $\text{Re}(z) \longrightarrow -\infty$ is ensured by the fact that the degree of the denominator $D(z)$ exceeds that of the numerator $N(z)$. Imposing the additional requirement that the rational function be a fourth-order approximation to the exponential leads to the coefficients listed in Table~\ref{tab_etd4_para}. The weights $w_i$ correspond to the partial fraction representation of $R(z)$ in Eq.~\ref{eqn2.32}. The key advantage of the RDP rational function over L-acceptable Padé approximations of comparable order is its real and distinct poles, which eliminate the need for complex arithmetic in linear systems.

\begin{table}
\begin{center}
\begin{tabular}{ c c c c} 
 \hline
  $k$ & $b_k$ & $w_k$& $a_k$ \\
  \hline
  \hline
 1 & 0.4751834017787114 & 20.10707940496431 & 1.579627631895178\\ 
 2 & 1.0000000000000000 & 0.5229558818011362 & 0.303085087930133\\ 
 3 & 0.3888888888888889 & -15.21083750434353 & -0.324250474367700\\
 4 & 0.7155553412275962 & -4.419197782421921 & \\
 \hline
\end{tabular}
\caption{Parameter values for the Partial fraction form of the RDP rational function.}\label{tab_etd4_para}
\end{center}
\end{table}


The corresponding matrix notation of the RDP function for the matrix exponential is given as 
\begin{equation}\label{eqn2.33}
R(-kA) = (I - a_1 kA + a_2 k^2A^2 - a_3 k^3 A^3) \left( \prod_{i=1}^4 (I + b_i k A) \right)^{-1},    
\end{equation}
with corresponding partial fraction decomposition,
\begin{equation}\label{eqn2.34}
  R(-kA) = \sum_{i=1}^4 (w_i)(I + b_i k A)^{-1}. 
\end{equation}
\noindent
The RDP approximation to $\mathrm{e}^{-\frac{k}{2}A}$ in matrix notation is similarly given by 
\begin{equation}\label{eqn2.35}
    R(-kA/2) = \left(I - \frac{k}{2} a_1A + \frac{k^2}{4}a_2A^2 - \frac{k^3}{8}a_3A^3\right)\left( \prod_{i=1}^4 \left(I + \frac{b_i}{2} k A \right) \right)^{-1}
\end{equation}  
\noindent
with partial fraction decomposition, 
\begin{equation}\label{eqn2.36}
    R (-kA/2) = \sum_{i = 1}^{4} w_i\left( I + \frac{b_i}{2}kA \right)^{-1}.
\end{equation}

\subsection{Derivation of the ETDRK4RDP-IF scheme}\label{etdrk4rdpif}
To begin, define the partial fraction form of the RDP rational approximation of the system matrices $A_1$ and $A_2$ using Eq.~\eqref{eqn2.34} as 

\begin{eqnarray}\label{eqn2.37}
  R(-kA_1) &=& \sum_{i=1}^4 (w_i)(I + b_i k A_1)^{-1} \\ 
  R(-kA_2) &=& \sum_{i=1}^4 (w_i)(I + b_i k A_2)^{-1}.
\end{eqnarray}

\noindent
We begin the derivation by considering each of the intermediate solutions $a_n, b_n, c_n$ and the main solution $U_{n+1}$ separately and simplify their results when the matrix exponentials are approximated using the respective RDP rational function. First, consider $\tilde{P}(kA)$ as defined in Eq.~\eqref{eqn2.2T} and approximate the matrix exponential $\mathrm{e}^{-kA/2}$ with the RDP function $R(-kA/2)$ in Eq.~\eqref{eqn2.35}. Thus, 
$$  \tilde{P}(kA) = - A^{-1}(\mathrm{e}^{-\frac{k}{2}A}-I). $$
So, 

\begin{eqnarray*}
\tilde{P}(kA) &=& - A^{-1}(R(-kA/2)- I)\\
\tilde{P}(kA) &=& -A^{-1}\left( (I - \frac{k}{2} a_1A + \frac{k^2}{4}a_2A^2 - \frac{k^3}{8}a_3A^3)\left( \prod_{i=1}^4 \left(I + \frac{b_i}{2} k A\right) \right)^{-1} - I \right). 
\end{eqnarray*}
\noindent
Then, factor out the matrix inverses, simplify algebraically, and write the results as a partial fraction decomposition, 

\begin{equation}\label{eqn2.44421}
     \tilde{P}(kA) = \sum_{i = 1}^{4} q_i\left( I + \frac{b_i}{2}kA \right)^{-1},
\end{equation}

where 
$$\displaystyle q_i = \frac{k \left( \mu_1 - 2\frac{\mu_2}{b_i} + 4\frac{\mu_3}{b_i^2} - 8\frac{\mu_4}{b_i^3} \right)}{\prod_{j=1, j\ne i}^{4} \left( 1 - \frac{b_j}{b_i}\right) } \quad \quad i = 1,2,3,4, $$

and, 
\begin{eqnarray*}
    \mu_1 &=& \frac{1}{2} (a_1 - (b_1 +b_2+b_3+b_4)),\ \mu_2 = \frac{1}{4} (a_2 - (b_1b_2 + b_3b_4 + (b_1 + b_2)(b_3+b_4)) \\
    \mu_3 &=& \frac{1}{8} (a_3 - ((b_1 + b_2)b_3b_4 + (b_3 + b_4)b_1b_2)),\ \mu_4 = -\frac{1}{16}(b_1b_2b_3b_4). 
\end{eqnarray*}


Consider the expression for $a_n$ in the ETDRK-IF scheme in Eq.~\eqref{eqn2.2} and approximate the matrix exponentials $\mathrm{e}^{-\frac{k}{2}A_1}$ and $ \mathrm{e}^{-\frac{k}{2}A_2}$ with their respective RDP rational approximation, $R(-kA_1/2)$ and $R(-kA_2/2)$. That is, 

$$a_n = \mathrm{e}^{-\frac{k}{2}A_2} \mathrm{e}^{-\frac{k}{2}A_1} U_n +\tilde{P}(kA_2) \mathrm{e}^{-\frac{k}{2}A_1}  F(U_n,t_n).$$
\noindent
By the commutativity of $\mathrm{e}^{-\frac{k}{2}A_2}$, $\mathrm{e}^{-\frac{k}{2}A_1}$ and $\tilde{P}(kA_2)$, we have that  
\begin{eqnarray*}
     a_n &=& \mathrm{e}^{-\frac{k}{2}A_1}[\mathrm{e}^{-\frac{k}{2}A_2}  U_n +\tilde{P}(kA_2)F(U_n,t_n)],\\
     a_n &=&  R(-kA_1/2)[  R(-kA_2/2) U_n +\tilde{P}(kA_2)F(U_n,t_n)],\\
     a_n &=&  \sum_{i = 1}^{4} w_i\left( I + \frac{b_i}{2}kA_1 \right)^{-1} (\tilde{a}_n),
\end{eqnarray*}
\noindent
where, 
$$ \tilde{a}_n =  R(-kA_2/2) U_n +\tilde{P}(kA_2)F(U_n,t_n).$$
By introducing the partial fraction form of the rationals, we have that 
\begin{equation*}
 \tilde{a}_n = \sum_{i = 1}^{4}  J_{ia} \quad \text{where}, \quad J_{ia} = \left(  I + \frac{b_i}{2}kA_2 \right)^{-1} (w_iU_n - q_iF(u_n, t_n)), \ i = 1,2,3,4.    
\end{equation*}
\noindent
Hence, to find $a_n$, we compute 
\begin{equation}
     a_n = \sum_{i = 1}^{4}  \tilde{J}_{ia},
\end{equation}
\noindent
where
$$ \tilde{J}_{ia} = w_i\left( I + \frac{b_i}{2}kA_1 \right)^{-1} (\tilde{a}_n), \quad i = 1,2,3,4.  $$


Similarly, consider the intermediate solution $b_n$ and approximate the matrix exponentials $\mathrm{e}^{-\frac{k}{2}A_1}$ and $ \mathrm{e}^{-\frac{k}{2}A_2}$ with their respective RDP rational approximations, $R(-kA_1/2)$ and $R(-kA_2/2)$. That is, 

\begin{eqnarray*}
    b_n &=& \mathrm{e}^{-\frac{k}{2}A_2} \mathrm{e}^{-\frac{k}{2}A_1} U_n +\tilde{P}(kA_2)F \left(a_n,t_n +\frac{k}{2} \right),\\
    b_n &=& \sum_{i = 1}^{4} w_i\left( I + \frac{b_i}{2}kA_2 \right)^{-1} \sum_{i = 1}^{4} w_i\left( I + \frac{b_i}{2}kA_1 \right)^{-1} U_n\\
     && - \sum_{i = 1}^{4} q_i\left( I + \frac{b_i}{2}kA_2 \right)^{-1} F \left(a_n,t_n +\frac{k}{2} \right).
\end{eqnarray*}
Set $\displaystyle J_{bn} =  \sum_{i = 1}^{4} w_i\left( I + \frac{b_i}{2}kA_1 \right)^{-1} U_n $ and compute its solution as 

$$ J_{bn} =  \sum_{i = 1}^{4} \tilde{J}_{bni} \quad \text{where} \quad \tilde{J}_{bni} = \left( I + \frac{b_i}{2}kA_1 \right)^{-1} w_i U_n.  $$
Thus, 
\begin{eqnarray*}
     b_n &=& \sum_{i = 1}^{4} w_i\left( I + \frac{b_i}{2}kA_2 \right)^{-1}[J_{bn}] - \sum_{i = 1}^{4} q_i\left( I + \frac{b_i}{2}kA_2 \right)^{-1} F \left(a_n,t_n +\frac{k}{2} \right), \\
      b_n &=& \sum_{i = 1}^{4} \left( I + \frac{b_i}{2}kA_2 \right)^{-1} \left[ w_i J_{bn} - q_iF \left(a_n,t_n +\frac{k}{2} \right) \right].\\
\end{eqnarray*}
Thus, to solve for $b_n$, we set 
$$\displaystyle b_n = \sum_{i = 1}^{4} \tilde{J}_{ib}, \quad \text{where}$$
$$ \tilde{J}_{ib} =  \left( I + \frac{b_i}{2}kA_2 \right)^{-1} \left[ w_i J_{bn} - q_iF \left(a_n,t_n +\frac{k}{2} \right) \right], \quad i = 1,2,3,4. $$

To evaluate $c_n$, we approximate the matrix exponentials $\mathrm{e}^{-\frac{k}{2}A_1}, \mathrm{e}^{-\frac{k}{2}A_2}$ and $\mathrm{e}^{-kA_1}$ with the corresponding RDP rational approximation $R(-kA_1/2), R(-kA_2/2)$ and $R(-kA_1)$. That is, 

\footnotesize
\begin{eqnarray*}
     c_n &=& \mathrm{e}^{-\frac{k}{2}A_2} \mathrm{e}^{-\frac{k}{2}A_1} a_n +\tilde{P}(kA_2) \left[ 2 \mathrm{e}^{-\frac{k}{2}A_1} F\left(b_n,t_n+\frac{k}{2} \right)- \mathrm{e}^{-\frac{k}{2}A_1}F(U_n,t_n) \right], \nonumber\\
    c_n &=& \mathrm{e}^{-\frac{k}{2}A_1} \left[\mathrm{e}^{-\frac{k}{2}A_2} a_n + 2\tilde{P}(kA_2) F \left(b_n,t_n+\frac{k}{2} \right) \right] - \tilde{P}(kA_2) \mathrm{e}^{-kA_1}F(U_n,t_n), \nonumber\\
    c_n &=&  \sum_{i = 1}^{4} w_i\left( I + \frac{b_i}{2}kA_1 \right)^{-1} \left[  \sum_{i = 1}^{4} w_i\left( I + \frac{b_i}{2}kA_2 \right)^{-1} a_n + 2 \sum_{i = 1}^{4} q_i\left( I + \frac{b_i}{2}kA_1 \right)^{-1} F \left(b_n,t_n+\frac{k}{2} \right) \right]\\
    && -  \sum_{i = 1}^{4} q_i\left( I + \frac{b_i}{2}kA_1 \right)^{-1}  \sum_{i = 1}^{4} w_i\left( I + b_ikA_1 \right)^{-1}F(U_n,t_n), \nonumber\\
    c_n &=&  \sum_{i = 1}^{4} w_i\left( I + \frac{b_i}{2}kA_1 \right)^{-1} [J_{cn1}] - \sum_{i = 1}^{4} q_i\left( I + \frac{b_i}{2}kA_2 \right)^{-1}  [J_{cn2}] \nonumber\\
\end{eqnarray*}

\normalsize
\noindent
where,
$$ J_{cn1} = \sum_{i = 1}^{4} w_i\left( I + \frac{b_i}{2}kA_2 \right)^{-1} a_n + 2 \sum_{i = 1}^{4} q_i\left( I + \frac{b_i}{2}kA_1 \right)^{-1} F \left(b_n,t_n+\frac{k}{2} \right),$$  
and 
$$J_{cn2} = \sum_{i = 1}^{4} w_i\left( I + b_ikA_1 \right)^{-1}F(U_n,t_n).$$

\noindent
$ J_{cn1}$ and $ J_{cn2}$ are evaluated as follows, 
$$ J_{cn1} =  \sum_{i = 1}^{4} \tilde{J}_{cn1_i} \quad \text{where}, \quad \tilde{J}_{cn1_i} = \left( I + \frac{b_i}{2}kA_2 \right)^{-1} (w_ia_n - 2q_iF(b_n, t_n + k/2)),$$
and 
$$ J_{cn2} =  \sum_{i = 1}^{4} \tilde{J}_{cn2_i}\quad \text{where}, \quad  \tilde{J}_{cn2_i} = \left( I + b_ikA_1 \right)^{-1} w_i F(U_n, t_n), \ i = 1,2,3,4.$$

\noindent
To simplify further, we set 
$$ J_{cn3} =  \sum_{i = 1}^{4} w_i\left( I + \frac{b_i}{2}kA_1 \right)^{-1} (J_{cn1}) \quad \text{and} \quad J_{cn4} = \sum_{i = 1}^{4} q_i\left( I + \frac{b_i}{2}kA_2 \right)^{-1}  (J_{cn2}). $$
Then $\displaystyle c_n = J_{cn3} + J_{cn4}$. We compute these intermediate approximations as follows, 
$$ J_{cn3} =  \sum_{i = 1}^{4} \tilde{J}_{cn3_i} \quad \text{where}, \quad  \tilde{J}_{cn3_i} = \left( I + \frac{b_i}{2}kA_1 \right)^{-1} w_i [J_{cn1}]$$
and 
$$ J_{cn4} =  \sum_{i = 1}^{4} \tilde{J}_{cn4_i} \quad \text{where}, \quad  \tilde{J}_{cn4_i} = \left( I + \frac{b_i}{2}kA_2 \right)^{-1} w_i [J_{cn2}]$$


To solve for the main solution, consider $U_{n+1}$ in Eq.~\eqref{eqn2.2} and simplify as follows
\begin{eqnarray*}
  U_{n+1} &=& \mathrm{e}^{-\frac{k}{2}A_1} \mathrm{e}^{-kA_2} U_n + P_1(kA_2) \mathrm{e}^{-\frac{k}{2}A_1}F(U_n,t_n), \nonumber \\
          && + 2P_2(kA_2) \mathrm{e}^{-\frac{k}{2}A_1} G \left( a_n, b_n, t_n + \frac{k}{2} \right) + P_3(kA_2) F(c_n,t_n+k), \\
  U_{n+1} &=& \sum_{i=1}^4 w_i \left( I + b_ikA_1 \right)^{-1} \left[ \sum_{i=1}^4 w_i\left( I + b_ikA_2 \right)^{-1} U_n + \sum_{i=1}^4 r_i\left( I + b_ikA_2 \right)^{-1} F(U_n,t_n)\right]\\
          && + \ 2\sum_{i=1}^4 g_i\left( I + b_ikA_2 \right)^{-1} \sum_{i=1}^4 w_i\left( I + \frac{b_i}{2}kA_1 \right)^{-1}G \left( a_n, b_n, t_n + \frac{k}{2} \right)  \\
          && + \sum_{i=1}^4 h_i\left( I + b_ikA_2 \right)^{-1} F(c_n,t_n+k).     
\end{eqnarray*}
\noindent
Set $$ U_{n1} = \sum_{i=1}^4 \left( I + b_ikA_2 \right)^{-1} [w_iU_n + r_iF(U_n, t_n)], $$
and $$ U_{n2} = 2 \sum_{i=1}^4 w_i\left( I + \frac{b_i}{2}kA_1 \right)^{-1} G \left( a_n, b_n, t_n + \frac{k}{2} \right).$$
Then,
\begin{equation}\label{eqn2.38}
     U_{n+1} = \sum_{i=1}^4 w_i \left( I + b_ikA_1 \right)^{-1} U_{n1} + \sum_{i=1}^4 g_i \left( I + b_ikA_2 \right)^{-1} [U_{n2} + h_i F(c_n,t_n+k)].
\end{equation}
The intermediate solutions $ U_{n1}$ and $ U_{n2}$ are obtained as follows,
\begin{eqnarray*}
    U_{n1} &=& \sum_{i=1}^4 \tilde{U}_{n1_i} \quad \text{where}, \ \tilde{U}_{n1_i} = \left( I + b_ikA_2 \right)^{-1} [w_iU_n + r_iF(U_n, t_n)], \quad i = 1,2,3,4.\\
     U_{n2} &=& \sum_{i=1}^4 \tilde{U}_{n2_i} \quad \text{where}, \ \tilde{U}_{n2_i} = 2\left( I + \frac{b_i}{2}kA_1 \right)^{-1}g_iw_i G \left( a_n, b_n, t_n + \frac{k}{2} \right).
\end{eqnarray*}
\noindent
By letting 
$$\displaystyle U_{n3} = \sum_{i=1}^4 w_i \left( I + b_ikA_1 \right)^{-1} U_{n1},$$
and
$$\displaystyle U_{n4} = \sum_{i=1}^4 g_i \left( I + b_ikA_2 \right)^{-1} [U_{n2} + h_i F(c_n,t_n+k)], $$
in Eq.~\eqref{eqn2.38}, we solve for $U_{n+1}$ as
$$ U_{n+1} = U_{n3} + U_{n4} $$
where $U_{n3}$ and $U_{n4}$  are computed as follows 
\begin{eqnarray*}
    U_{n3} &=& \sum_{i=1}^4 \tilde{U}_{n3_i} \quad \text{where}, \ \tilde{U}_{n3_i} = \left( I + b_ikA_1 \right)^{-1}U_{n1}, \quad i = 1,2,3,4.\\
     U_{n4} &=& \sum_{i=1}^4 \tilde{U}_{n4_i} \quad \text{where}, \ \tilde{U}_{n4_i} =\left( I + b_ikA_2 \right)^{-1} [U_{n2} + h_i F(c_n,t_n+k)], \ i = 1,2,3,4.
\end{eqnarray*}

\subsection{Implementation of the ETDRK4RDP-IF Algorithm}
\noindent
The ETDRK4RDP-IF scheme has 48 linear solves at each time step. The intermediate solutions $a_n$ and $b_n$ have 8 linear solves each, while $c_n$ and the main solution $U_{n+1}$ have 16 linear solves. The fully discrete ETDRK4RDP-IF scheme is implemented in serial with the following steps:
\begin{enumerate}
    \item (a) Solve the system  
    $$\displaystyle \left( I + \frac{b_i}{2}kA_2 \right)J_{ia} = (w_iU_n - q_iF(U_n, t_n)) \quad \text{for} \quad i = 1,2,3,4.$$\label{A11} \\
    (b) Set $\displaystyle \tilde{a}_n = \sum_{i=1}^4 J_{ia}$.\label{A12}\\
    (c) Solve the system $\displaystyle \left( I + \frac{b_i}{2}kA_1 \right)\tilde{J}_{ia} = w_i \tilde{a_n}$ for $i = 1,2,3,4$.\label{A13}\\
    (d) Set $\displaystyle a_n = \sum_{i=1}^4 \tilde{J}_{ia}$.\label{A14}

    \item (a) Solve the system $\displaystyle \left( I + \frac{b_i}{2}kA_1 \right) \tilde{J}_{bni} = w_iU_n $, for $i = 1,2,3,4$.\label{B11}\\
     (b) Set $\displaystyle J_{bn} = \sum_{i=1}^4 \tilde{J}_{bni}$.\label{B12}\\
     (c) Solve the system $\displaystyle \left( I + \frac{b_i}{2}kA_2 \right) \tilde{J}_{ib} = \left[ w_i J_{bn} - q_iF \left(a_n,t_n +\frac{k}{2} \right) \right] $ for $i = 1,2,3,4$.\label{B13}\\
     (d) Set $\displaystyle b_n = \sum_{i=1}^4 \tilde{J}_{ib}$.\label{B14}
     
    \item (a) Solve the system 
    $$\displaystyle \left( I + \frac{b_i}{2}kA_2 \right) \tilde{J}_{cn1_i} = w_ia_n - 2q_i[2F(b_n, t_n + k/2)] \quad \text{for} \quad i = 1,2,3,4$$.\label{C11}\\
    (b) Set $\displaystyle J_{cn1} = \sum_{i=1}^4 \tilde{J}_{cn1_i} $.\label{C12}\\
    (c) Solve the system $\displaystyle \left( I + b_ikA_1 \right) \tilde{J}_{cn2_i} = w_iF(U_n, t_n)$ for $i = 1,2,3,4$.\label{C13}\\
    (d) Set $\displaystyle J_{cn2} = \sum_{i=1}^4 \tilde{J}_{cn2_i} $.\label{C14}\\
    (e) Solve the system $\displaystyle \left( I + \frac{b_i}{2}kA_1 \right) \tilde{J}_{cn3_i} = w_i J_{cn1}$ for $i = 1,2,3,4$.\label{C15}\\
    (f) Set $\displaystyle J_{cn3} = \sum_{i=1}^4 \tilde{J}_{cn3_i} $.\label{C16}\\
    (g) Solve the system $\displaystyle \left( I + \frac{b_i}{2}kA_2 \right) \tilde{J}_{cn4_i} = q_i J_{cn2}$ for $i = 1,2,3,4$.\label{C17}\\
    (h) Set $\displaystyle J_{cn4} = \sum_{i=1}^4 \tilde{J}_{cn4_i} $.\label{C18}\\
    (i) Compute $c_n$ as $$\displaystyle c_n =  J_{cn3} + J_{cn4}.$$

    \item (a) For the main solution, solve the system 
    $$ ( I + b_ikA_2 ) \tilde{U}_{n1_i} = w_iU_n + r_iF(U_n, t_n) \quad \text{for} \ i = 1,2,3,4.$$\label{u11}
    \noindent
    (b) Set $\displaystyle U_{n1} = \sum_{i=1}^4 \tilde{U}_{n1_i} $.\label{u12}\\
    (c) For $i = 1,2,3,4$, solve $\displaystyle \left( I + \frac{b_i}{2}kA_2 \right) \tilde{U}_{n2_i} = 2w_i G \left( a_n, b_n, t_n + \frac{k}{2} \right)$.\label{u13}\\
    (d) Set $\displaystyle U_{n2} = \sum_{i=1}^4 \tilde{U}_{n2_i} $.\label{u14}\\
    (e) Solve the system $\displaystyle \left( I + b_ikA_1 \right) \tilde{U}_{n3_i} = w_i U_{n1} $ for $i = 1,2,3,4$.\label{u15}\\
    (f) Set $\displaystyle U_{n3} = \sum_{i=1}^4 \tilde{U}_{n3_i} $.\label{u16}\\
    (g) Solve the system 
    $$\displaystyle \left( I + b_ikA_2 \right) \tilde{U}_{n4_i} =  U_{n2} + h_i F(c_n, t_n + h) \quad \text{for} \ i = 1,2,3,4.$$\label{u17}\\
    (h) Set $\displaystyle U_{n4} = \sum_{i=1}^4 \tilde{U}_{n4_i} $.\label{u18}\\
    (g) Then set,  $$U_{n+1} = U_{n3} + U_{n4}$$.
\end{enumerate}

\noindent
The algorithm can be implemented in a parallel environment using four processors to solve the linear systems associated with each stage concurrently. We investigate the performance of the new scheme by comparing it with the standard fourth-order ETDRK4RDP and ETDRK4P22-IF schemes.

\section{Numerical Experiments and Results}\label{section3}

To validate the accuracy, convergence, and computational efficiency of the proposed ETDRK4RDP-IF scheme, we performed a variety of numerical experiments on a representative set of benchmark problems with varying initial data and boundary conditions. The new method was compared with two existing fourth-order ETD schemes, namely ETDRK4P22-IF (with dimensional splitting) and ETDRK4RDP (without dimensional splitting), introduced in \cite{asante2025fourth} and \cite{attipoe2025fourth} respectively. The ETDRK4P22-IF scheme employs a Padé$(2,2)$ rational approximation for the matrix exponential, which is not L-acceptable and shows limited damping of spurious oscillations. In contrast, the ETDRK4RDP method utilizes an RDP approximation within the ETDRK framework of Cox and Matthews \cite{cox2002exponential}, but is not dimensionally split. In this way, we are able to evaluate the value of dimensional splitting as well as the RDP rational approximation. The fully discrete ETDRK4RDP scheme (see \cite{attipoe2025fourth}) is implemented in parallel with four processors. The fully discrete ETDRK4P22-IF scheme (see \cite{asante2025fourth}) is implemented in parallel with up to three processors. 

The numerical experiments were conducted in two spatial dimensions and encompass both scalar and coupled systems of reaction–diffusion equations subject to homogeneous Dirichlet and Neumann boundary conditions, thereby providing an evaluation of the proposed scheme across a broad class of problems relevant to computational physics. The spatial discretization was done using a fourth order finite difference scheme, as discussed in \cite{attipoe2025fourth}. The error, $E(k)$ at time step $k$, is measured with the $L_{\infty}-$ norm, and the estimated order of convergence calculated using the formula     
$$ p = \frac{\log[E(k)/E(k/2)]}{\log 2}. $$
\noindent
Calculation of the error for problems with and without exact solution was done in a similar way as was discussed in \cite{attipoe2025fourth}. 

Serial implementation of all algorithms for the numerical experiments was run in MATLAB R2023b on a Dell Latitude 5440 with 13th Gen Intel(R) Core(TM) i7-1355U 1.70 GHz with 32 GB RAM. The parallel variants of the algorithms were in Fortran using OpenMP with 2-4 processors on a PC with 32 $\times$ 12th Gen Intel Core i9-14900 with 32 GiB of RAM and 1.0TB of memory running Ubuntu 22.04. The resulting sparse linear systems were solved using the UMFPACK library \cite{Davis2025suitesparse}, with matrices stored in compressed column storage format \cite{Hanyk2014mUMFPACK}. Parallelization was carried out in Fortran instead of MATLAB due to MATLAB's reliance on multithreading which limits the efficiency of explicit parallelization. The MATLAB and Fortran implementations of the ETDRK4RDP-IF scheme for all numerical experiments presented in this work are publicly available at \url{https://github.com/wisdomattipoe/ETDRK4RDP-IF-Scheme-Codes.git}

\subsection{Model problem with a homogeneous Dirichlet boundary condition}\label{Dirichlet}
\begin{figure}[h]
\centering
\includegraphics[width=12cm]{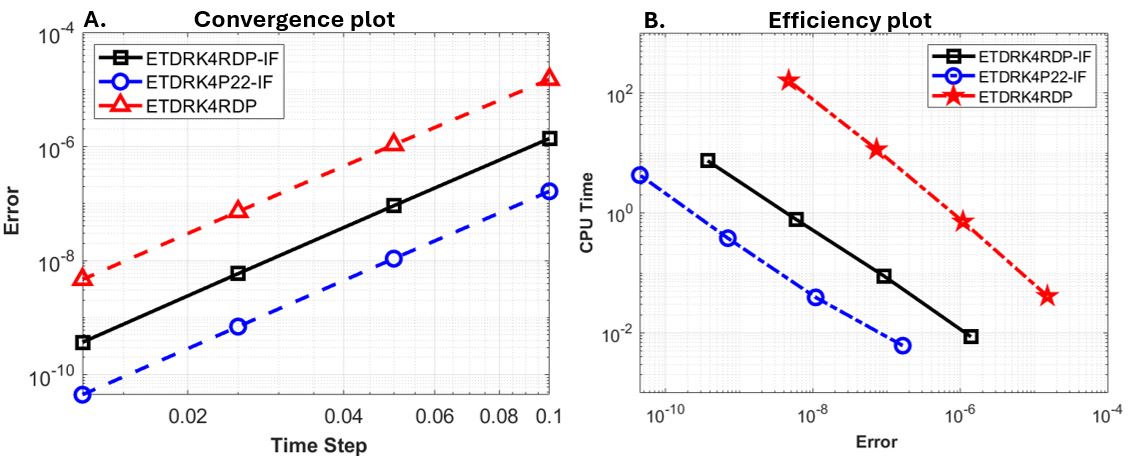}
\caption{\text{Convergence and efficiency plots for test problem \ref{Dirichlet}}. (A) Convergence and (B) Efficiency of the fourth-order ETDRK4RDP-IF scheme compared with the existing fourth-order ETDRK4P22-IF, and ETDRK4RDP schemes.} 
\label{Fig. 3.1}
\end{figure}

\noindent
The performance of the proposed method is examined using the following two-dimensional reaction–diffusion test-system, equipped with homogeneous Dirichlet boundary conditions. The system has a closed-form solution given by $\displaystyle u(x, y, t) = \mathrm{e}^{-3t}\cos(x)\cos(y) $. 

\begin{eqnarray}
\frac{\partial u}{\partial t}  &=& \Delta u - u, \quad \frac{-\pi}{2} < x, y < \frac{\pi}{2}, \quad t \in [0, 1], \\
u(x, y, 0) &=& \cos(x) \cos(y) \nonumber.
\end{eqnarray}
\noindent
Spatial discretization of the negative Laplacian is performed as discussed in Chapter 2 of \cite{asante2025fourth}. The convergence plot in Figure~\ref{Fig. 3.1}A and grid refinement study in  Table~\ref{table3.1} confirm the fourth-order accuracy of the proposed scheme. For all time steps examined, the serial implementation of the ETDRK4RDP-IF scheme in MATLAB is more accurate than the existing ETDRK4RDP scheme but less accurate than the ETDRK4P22-IF scheme. See Figure~\ref{Fig. 3.1} and Table~\ref{table3.1}.

Remarkably, although the ETDRK4RDP-IF scheme performs four times the number of linear solves as the standard ETDRK4RDP scheme it is about 20 times faster at the finest resolutions with $k = 0.0125$ and $h = 0.01$ (See Table~\ref{table3.1}). Comparing the speed of serial implementations of the two dimensionally split methods, ETDRK4RDP-IF and ETDRK4P22-IF, we see that ETDRK4P22-IF is two times faster than our proposed scheme. This makes sense as ETDRK4RDP-IF performs twice as many linear solves as the existing ETDRK4P22-IF scheme. At the finest spatial resolution, the parallel implementations of both schemes, ETDRK4RDP-IF and ETDRK4P22-IF, exhibit comparable computational performance.

For all time steps considered, the standard ETDRK4RDP scheme exhibits lower accuracy and incurs significantly higher CPU time relative to ETDRK4RDP-IF (See Table~\ref{table3.1}). The increased computational cost is primarily attributable to the solution of larger linear systems arising within the formulation of the method. A substantial performance enhancement is observed for the proposed ETDRK4RDP-IF method when its parallel implementation is compared with the standard parallel ETDRK4RDP scheme at the finest resolution, achieving a speedup exceeding a factor of seven (see, Table~\ref{table3.1}).  

\begin{table}[ht]
 \begin{center}
 \resizebox{\textwidth}{!}{
\begin{tabular}{ p{3.2cm}|p{1cm}|p{1cm}|p{2.2cm}|p{1cm}|p{2cm}|p{1.2cm}|p{2cm} }
 \hline
 \hline
   Method    & $k$    & $h$  & Error                 & $p$  & CPU Time MATLAB& 1-core& Multicore\\
   \hline
   \hline
ETDRK4RDP-IF & 0.10   & 0.08 & $1.38 \times 10^{-6}$ & $-$  & 0.035           & 0.018  & 0.005  \\       
             & 0.05   & 0.04 & $9.15 \times 10^{-8}$ & 3.92 & 0.118           & 0.163  & 0.041  \\
             & 0.025  & 0.02 & $5.89 \times 10^{-9}$ & 3.96 & 0.918           & 1.468  & 0.367  \\   
             & 0.0125 & 0.01 & $3.70 \times 10^{-10}$& 3.96 & 9.343           & 14.279 & 3.570  \\
 \hline
 \hline
   ETDRK4RDP & 0.10   & 0.08 & $1.50 \times 10^{-5}$ & $-$  & 0.062           & 0.031  & 0.010  \\       
             & 0.05   & 0.04 & $1.07 \times 10^{-6}$ & 3.80 & 1.037           & 0.381  & 0.143  \\
             & 0.025  & 0.02 & $7.23 \times 10^{-8}$ & 3.89 & 17.152          & 4.842  & 2.391  \\   
             & 0.0125 & 0.01 & $4.66 \times 10^{-9}$ & 3.96 & 178.107         & 53.832 & 27.501 \\ 
 \hline
 \hline            
ETDRK4P22-IF & 0.10   & 0.08 & $1.64 \times 10^{-7}$ & $-$  & 0.021           & $0.013$  & $0.006$ \\      
             & 0.05   & 0.04 & $1.08 \times 10^{-8}$ & 3.92 & 0.071           & $0.114$  & $0.025$ \\
             & 0.025  & 0.02 & $6.96 \times 10^{-10}$& 3.96 & 0.376           & $1.115$  & $0.179$ \\  
             & 0.0125 & 0.01 & $4.46 \times 10^{-11}$& 3.96 & 4.598           & $6.916$  & $3.360$ \\ 
 \hline
 \hline      
\end{tabular}
  }
  \caption{\label{table3.1} Convergence table showing performance of the new ETDRK4RDP-IF scheme compared with the standard ETDRK4RDP and ETDRK4P22-IF schemes for the model test problem in Example \ref{Dirichlet} with $T =1$. The table also displays the CPU time for serial implementation of all schemes in MATLAB, as well as single-core (1-core) and multicore implementation in Fortran.}
\end{center}
\end{table}


\subsection{Nonlinear enzyme kinetics with mismatched initial and boundary data}\label{mismatched}
Next, we evaluate the performance of the ETDRK4RDP-IF scheme on the well-known non-smooth enzyme kinetics model subject to homogeneous Dirichlet conditions. The model constitutes a two-dimensional extension of the one-dimensional reaction–diffusion equation incorporating Michaelis–Menten enzyme kinetics \cite{muller2015methods, bhatt2015locally}. This formulation represents the prototypical enzyme–substrate interaction, describing an enzyme-catalyzed conversion of a single substrate into a single product \cite{muller2015methods}. The model is presented mathematically as,

$$ \frac{\partial u}{\partial t} = d \left(\frac{\partial^2u}{\partial x^2} +\frac{\partial^2u}{\partial y^2} \right) - \frac{u}{(1+u)}\,,\quad 0<x\,,\,y<1\,,\quad t>0, $$
\noindent
where the diffusion coefficient $d = 1$ and has the following initial condition
$$u(x, y, 0) = 1.$$

\begin{figure}[h]
\centering
\includegraphics[width=12cm]{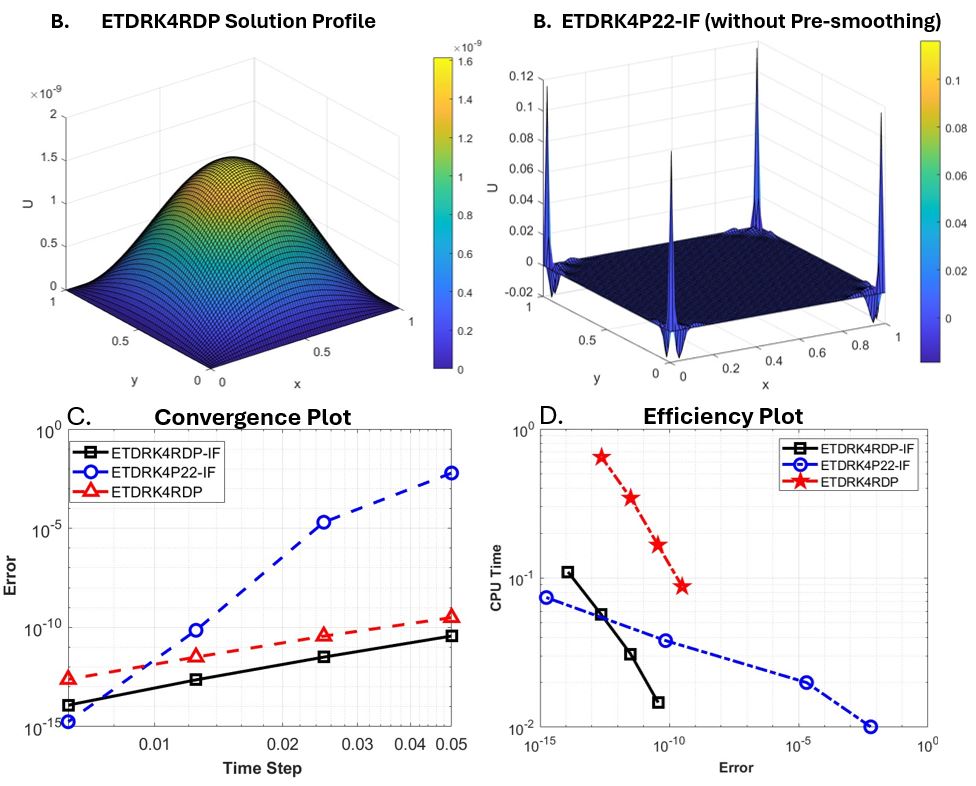}
\caption{\textbf{Convergence and efficiency plots for Example \ref{mismatched}}. (A) Numerical solution using ETDRK4RDP-IF (B) ETDRK4P22-IF solution without Pre-smoothing (C) Convergence Plot (D) Efficiency plot of the fourth-order ETDRK4RDP-IF scheme compared with the existing fourth-order ETDRK4P22-IF and the standard ETDRK4RDP schemes.} 
\label{Fig. 3.2}
\end{figure}

\noindent
Notably, the prescribed initial condition is discontinuous at the boundary where homogeneous Dirichlet conditions are imposed. Such incompatibilities between the initial and boundary data are well known to induce spurious oscillations which, if not adequately controlled during the early stages with appropriate pre-smoothing techniques, may substantially degrade the accuracy of the numerical solution \cite{asante2025fourth}. The system does not admit a known closed-form solution. Without pre-smoothing, the ETDRK4P22-IF scheme fails to damp out spurious oscillations for relatively large time steps (Figure~\ref{Fig. 3.2}B) and struggles to maintain fourth-order accuracy (Table~\ref{table3.2}). 
This phenomenon is also observed in \cite{asante2025fourth, attipoe2025fourth}. The dimensionally split ETDRK4RDP-IF is observed to be more efficient than the unsplit ETDRK4RDP in both serial and parallel implementations, gaining up to 5 times speed up in the multi-core environment. The L-acceptability of the ETDRK4RDP-IF scheme permits the smoothing of the boundary (Figure~\ref{Fig. 3.2} A), ensuring that the scheme converges with fourth-order accuracy (Table~\ref{table3.2}). In Figure~\ref{Fig. 3.2}C, we see the inconsistency in the convergence results of the ETDRK4P22-IF without pre-smoothing, whereas the L-acceptable schemes all show accurate performance. As can be observed from Table~\ref{table3.2} and Figure~\ref{Fig. 3.2}C, the accuracy of ETDRK4RDP-IF is comparable to that of ETDRK4RDP and performs about 6 times faster in CPU time when both schemes are implemented in serial. Parallelization results in a 3 times speed up for ETDRK4RDP-IF and a two times speed up for ETDRK4RDP, compared to their serial implementations in Fortran.

\begin{table}[ht]
 \begin{center}
 \resizebox{\textwidth}{!}{
\begin{tabular}{ p{3.2cm}|p{1.4cm}|p{1cm}|p{2.2cm}|p{1cm}|p{2cm}|p{1.2cm}|p{2cm} }
 \hline
 \hline
   Method    & $k$    & $h$   & Error                 & $p$  & CPU Time MATLAB& 1-core& Multicore\\
   \hline
   \hline
ETDRK4RDP-IF & 0.05   & 0.025 & $3.67 \times 10^{-11}$ & $-$  & 0.0146           & 0.017  & 0.005  \\       
             & 0.025  & 0.025 & $3.17 \times 10^{-12}$ & 3.53 & 0.0309           & 0.034  & 0.009  \\
             & 0.0125 & 0.025 & $2.26 \times 10^{-13}$ & 3.81 & 0.0572           & 0.074  & 0.019  \\   
             & 0.00625& 0.025 & $1.20 \times 10^{-14}$ & 4.14 & 0.1099           & 0.165  & 0.036  \\
 \hline
 \hline
   ETDRK4RDP & 0.05   & 0.025 & $3.13 \times 10^{-10}$ & $-$  & 0.0877          & 0.063  & 0.016 \\       
             & 0.025  & 0.025 & $3.62 \times 10^{-11}$ & 3.11 & 0.1663          & 0.129  & 0.032 \\
             & 0.0125 & 0.025 & $3.19 \times 10^{-12}$ & 3.50 & 0.3435          & 0.272  & 0.069 \\   
             & 0.00625& 0.025 & $2.39 \times 10^{-13}$ & 3.74 & 0.6468          & 0.484  & 0.122 \\ 
 \hline
 \hline            
ETDRK4P22-IF & 0.05   & 0.025 & $6.10 \times 10^{-3}$ & $-$   & 0.0101           & $0.016$  & $0.009$ \\    
             & 0.025  & 0.025 & $2.01 \times 10^{-5}$ & 8.24  & 0.0199           & $0.028$  & $0.036$ \\
             & 0.0125 & 0.025 & $7.13 \times 10^{-11}$& 18.11 & 0.0381           & $0.036$  & $0.049$ \\   
             & 0.00625& 0.025 & $1.74 \times 10^{-15}$& 15.32 & 0.0741           & $0.083$  & $0.042$ \\ 
 \hline
 \hline      
\end{tabular}
  }
  \caption{\label{table3.2} Convergence table showing fourth-order convergence of the new ETDRK4RDP-IF scheme compared with the existing ETDRK4RDP and ETDRK4P22-IF schemes for the model test problem in Example \ref{mismatched} with $T =1$. The table also displays the CPU time for serial implementation of all schemes in MATLAB, as well as single-core (1-core) and multicore implementation in Fortran.}
\end{center}
\end{table}

\subsection{Gray-Scott complex dynamics RDE }\label{grayscott}
Finally, we consider the coupled system of two-dimensional nonlinear reaction–diffusion equations, commonly referred to as the Gray–Scott model, which serves as a prototypical framework for the study of autocatalytic chemical kinetics and pattern formation \cite{doelman1997pattern}. The model is characterized by quadratic nonlinear reaction terms coupled with diffusive transport, leading to a rich spectrum of dynamical behaviors, including pulse splitting and self-replicating structures. Depending on the choice of diffusion coefficients and feed–removal parameters, the system exhibits steady states, oscillatory dynamics, and complex spatiotemporal patterns, making it a challenging benchmark problem for the analysis and numerical approximation of nonlinear PDEs. The partial differential equation modeling this process is given by
\begin{eqnarray*}
    \frac{\partial u}{\partial t} + d_1\Delta u &=& -uv^2 + \gamma(1 - u)\\
    \frac{\partial v}{\partial t} + d_2\Delta v &=&  uv^2 + (\gamma + \beta)v
\end{eqnarray*}
with a homogeneous Neumann boundary and an initial condition given by 
\begin{eqnarray*}
    u_0(x, y, 0) = 1 - \frac{1}{2}\left(\sin^{100}(\pi x) \sin^{100}(\pi y)\right) \\
    v_0(x, y, 0) = \frac{1}{4}\left(\sin^{100}(\pi x) \sin^{100}(\pi y)\right) 
\end{eqnarray*}

where $d_1 = 0.001, d_2 = 0.0001, \gamma = 0.024$ and $\beta = 0.06$. 

\begin{figure}[h]
\centering
\includegraphics[width=12cm]{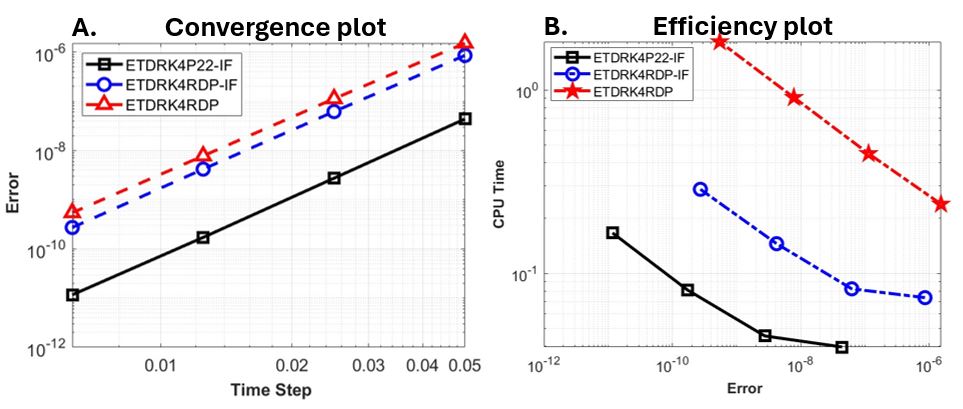}
\caption{\textbf{(A) Convergence and (B)  plots for Example~\ref{grayscott}} of the new fourth-order ETDRK4RDP-IF scheme compared with the existing fourth-order ETDRK4P22-IF and the standard ETDRK4RDP schemes.} 
\label{Fig. 3.5}
\end{figure}

\begin{table}[ht]
 \begin{center}
 \resizebox{\textwidth}{!}{
\begin{tabular}{ p{3.2cm}|p{1.2cm}|p{1.2cm}|p{2.2cm}|p{1cm}|p{2cm}|p{1.2cm}|p{2cm} }
 \hline
 \hline
   Method    & $k$     & $h$     & Error                & $p$  & CPU Time (Serial)& 1-core & Multi-core \\
   \hline
   \hline
ETDRK4RDP-IF & 0.05    & 0.025 & $8.59 \times 10^{-7}$  & $-$  & 0.074        & $0.083$  & $0.021$  \\                             & 0.025   & 0.025 & $6.21 \times 10^{-8}$  & 3.79 & 0.082        & $0.147$  & $0.037$  \\
             & 0.0125  & 0.025 & $4.20 \times 10^{-9}$  & 3.89 & 0.145        & $0.198$  & $0.075$  \\   
             & 0.00625 & 0.025 & $2.73 \times 10^{-10}$ & 3.94 & 0.288        & $0.363$  & $0.142$  \\
 \hline
 \hline
   ETDRK4RDP & 0.05    & 0.025 & $2.83 \times 10^{-6}$  & $-$  & 0.239       & $0.275$  & $0.069$ \\                            & 0.025   & 0.025 & $1.44 \times 10^{-7}$  & 3.75 & 0.450       & $0.511$  & $0.128$ \\
             & 0.0125  & 0.025 & $7.51 \times 10^{-9}$  & 3.86 & 0.910       & $0.935$  & $0.235$ \\   
             & 0.00625 & 0.025 & $4.13 \times 10^{-10}$ & 3.92 & 1.833       & $1.923$  & $0.482$ \\ 
 \hline
 \hline            
ETDRK4P22-IF & 0.05    & 0.025 & $8.34 \times 10^{-8}$  & $-$  & 0.040        & $0.049$  & $0.025$ \\               
             & 0.025   & 0.025 & $5.52 \times 10^{-9}$  & 4.00 & 0.046        & $0.086$  & $0.046$ \\
             & 0.0125  & 0.025 & $1.77 \times 10^{-10}$ & 4.00 & 0.081        & $0.168$  & $0.081$ \\   
             & 0.00625 & 0.025 & $1.72 \times 10^{-11}$ & 3.88 & 0.166        & $0.322$  & $0.168$ \\ 
 \hline
 \hline      
\end{tabular}
  }
  \caption{\label{table3.5} Convergence table showing fourth-order convergence of the new ETDRK4RDP-IF scheme compared with the existing ETDRK4RDP and ETDRK4P22-IF schemes for Example \ref{grayscott} with $T =1$. The table also displays the CPU time for serial implementation of all schemes in MATLAB, as well as single-core (1-core) and multicore implementation in Fortran.}
\end{center}
\end{table}

In the absence of an exact analytical solution, the reported errors correspond to coarse-to-fine grid errors, as defined in the introduction of this section. Empirical convergence results indicate that the proposed scheme achieves fourth-order temporal accuracy, exhibiting a level of accuracy comparable to that of the existing fourth-order ETDRK4RDP and ETDRK4P22-IF scheme (Table \ref{table3.5} and Fig. \ref{Fig. 3.5}). Notably, the proposed ETDRK4RDP-IF scheme demonstrates significantly improved computational efficiency relative to the standard ETDRK4RDP method. In particular, at the spatial mesh size $h = 0.0125$ and time step $k = 0.00625$, the ETDRK4RDP-IF scheme achieves approximately a six-fold reduction in CPU time under a serial implementation (Table \ref{table3.5} and Fig. \ref{Fig. 3.5}).
It is noteworthy that the ETDRK4RDP-IF scheme performs roughly five times as many linear solves per time step as the ETDRK4P22-IF scheme; nevertheless, the overall computational cost remains competitive when the serial CPU times are compared (Table \ref{table3.5}). As expected, all schemes considered exhibit comparable accuracy. A parallel implementation of the ETDRK4RDP-IF scheme using four cores yields a further reduction in computational cost, improving the CPU time by approximately a factor of 2 at the finest spatial and temporal resolutions (Table  \ref{table3.5}). Parallel implementation of the ETDRK4RDP-IF scheme still outperforms the ETDRK4RDP method by approximately 4 times in CPU time. This observation suggests that, for the Gray–Scott model considered, parallel implementations in Fortran provide the most favorable performance for these schemes.


\section{Conclusion}\label{section4}
In this work, we developed and tested a novel fourth-order, dimensionally split ETD Runge–Kutta scheme (ETDRK4RDP-IF) that combines the standard fourth-order finite differences for spatial discretizations with a real-and-distinct-pole (RDP) rational approximation for the corresponding matrix exponentials.

Benchmarked against the unsplit ETDRK4RDP and Padé(2,2)-based ETDRK4P22-IF scheme, the proposed method exhibits superior robustness for stiff reaction–diffusion equations with non-smooth initial conditions, owing to the $L$-acceptability of the RDP approximation. The serial implementation of the ETDRK4RDP-IF scheme shows a better computational efficiency than the standard ETDRK4RDP across all the test problems considered, even though it performs 4 times more linear solves. In some cases, we gain up to 18 times more speed-up at the finest resolutions. Parallelization of the new scheme also increases its computational efficiency by almost a factor of 3 in some cases. Across all test problems, we observed that parallelization in Fortran significantly improved the CPU time for ETDRK4RDP-IF, ETDRK4RDP and ETDRK4P22, with the most significant improvement observed for ETDRK4RDP-IF. Even though the ETDRK4RDP-IF scheme performs about 5 times more linear solves than ETDRK4P22-IF, they perform competitively on all the problems considered. For the RDP-type schemes, by distributing the four real linear systems that arise at each stage across four processors, the computational cost per processor is reduced to one real linear solve per stage, while avoiding complex arithmetic. In contrast, competing methods require solving a linear system involving complex-valued arithmetic. 

Future work will focus on extending these schemes to three-dimensional problems and incorporating adaptive time-stepping to enhance computational efficiency. On the theoretical side, we aim to establish rigorous proofs of fourth-order convergence and identify the requisite conditions on the differential operator and non-linear source term that guarantee both fourth-order accuracy and positivity preservation of the numerical solution.

\section{Acknowledgment}\label{section6}
The authors would like to thank the editor and referees whose suggestions led to a considerably improved version of the initial manuscript.

\section{Statements and Declarations}
\subsection*{Competing interests}
The authors declare that they have no known competing financial interests or personal relationships that could have appeared to influence the work reported in this paper.

\end{document}